\documentclass{pnastwo}

\usepackage[tbtags]{amsmath}
\usepackage{amssymb,amsfonts}
\usepackage[dvips]{graphicx}

\begin{document}



\conflictofinterest{Conflict of interest footnote placeholder}

\track{Insert 'This paper was submitted directly to the 
PNAS office.' when applicable.}




\title{Efficient Monte Carlo sampling by parallel marginalization}




\author{Jonathan Weare\thanks{E-mail: weare@math.berkeley.edu}
\affil{1}{Department of Mathematics, University of California, Berkeley  
CA 94720-3840 }}

\contributor{}

\maketitle

\begin{article}

\begin{abstract}Markov chain
Monte Carlo sampling methods often suffer from long correlation times.
Consequently, these methods must be run
for many steps to generate an independent sample.  In this paper a method is proposed to overcome this difficulty.
The method utilizes information from rapidly equilibrating coarse Markov chains 
that sample marginal distributions of the full system.  This is accomplished through
 exchanges between the full chain and the auxiliary coarse chains.
Results of numerical tests on 
the bridge sampling and filtering/smoothing problems for a stochastic differential
equation are presented.
\end{abstract}




\dropcap{I}n order to understand the behavior of a physical system it is often
necessary to generate samples from complicated high dimensional distributions.  
The usual tools for sampling from these distributions are Markov chain
Monte Carlo methods (MCMC) by which one constructs a Markov chain whose trajectory
averages converge to averages with respect to 
the distribution of interest.  For some simple systems it is 
possible to construct Markov chains with independent values at each
step.  In general, however,
spatial correlations in the system
of interest result in long correlation times in the Markov chain and hence
slow convergence of the chain's trajectory averages.
In this paper, a method is proposed to alleviate the difficulties caused by
spatial correlations in high dimensional systems.  
The method, parallel marginalization, is tested on
two stochastic differential equation conditional path sampling problems.

Parallel marginalization takes advantage
of the shorter correlation lengths present in
 marginal distributions of the target density.
Auxiliary 
Markov chains that sample approximate marginal distributions are evolved
simultaneously with the Markov chain that samples the distribution of interest.
By swapping their configurations, these auxiliary 
chains pass information between themselves and with the chain sampling
the original
distribution.  As shown below, these swaps are made
in a manner consistent with both the original distributions and the approximate 
marginal distributions.  The numerical examples indicate that improvement in efficiency of
parallel marginalization over standard MCMC techniques can be significant.

The design of efficient methods to approximate marginal 
distributions was addressed in \cite{c03} by Chorin and in \cite{s05}
 by Stinis.
The use of Monte Carlo updates on coarse subsets of variables is not
a new concept (see \cite{gs89} and the references therein).
The method presented in \cite{gs89} does not use marginal distributions.
However, attempts have been made previously to use 
marginal distributions to accelerate the
convergence of MCMC (see \cite{br01,o05}).
In contrast to parallel marginalization, the methods
proposed in \cite{br01} and \cite{o05} do not preserve the distribution
of the full system and therefore are not guaranteed to converge.
The parallel construction used here is motivated by the parallel tempering
method
(see \cite{l02}), and allows efficient comparison of the auxiliary chains and
the original chain.  See references \cite{l02} and \cite{bh02} for expositions
of standard MCMC methods.

Parallel marginalization for problems in Euclidean state spaces 
is described in detail in the next two sections.  In the final sections
the conditional path sampling problem is described and numerical results
are presented for the bridge sampling and smoothing/filtering problems.

\section{Parallel Marginalization}
For the purposes of the discussion in this section, we assume that appropriate 
approximate marginal distributions are available.
As discussed in a later section, they may be provided by
coarse models of the physical problem as in the examples below,
or they may be calculated via the methods in \cite{c03} and \cite{s05}.

Assume that
the $d_0$ dimensional system of interest has
a probability density, $\pi_0(x_0)$, where
$x_0\in\mathbb{R}^{d_0}$. 
Suppose further that, by the Metropolis-Hastings or any other method (see \cite{l02}),
we can construct a Markov chain, $Y_0^n\in\mathbb{R}^{d_0}$, which has
$\pi_0$ as its stationary measure.  That is, for two points $x_0,y_0\in\mathbb{R}^{d_0}$
$$
\int T_0(x_0\rightarrow y_0)\pi_0(x_0)\ dx_0 = \pi_0(y_0)
$$
where $T_0(x_0\rightarrow y_0)$ is the probability density of a move
to $\left\{Y_0^{n+1}=y_0\right\}$ given that $\left\{Y_0^n=x_0\right\}$.  
Here, $n$ is the algorithmic step.  Under
appropriate conditions (see \cite{l02}), 
averages over a trajectory of $\left\{Y_0^n\right\}$ will
converge to averages over $\pi_0$, i.e. for an objective function $g(x_0)$
$$
\frac{1}{N}\sum_{n=0}^{N-1} g\left(Y_0^n\right)
\rightarrow \mathbf{E}\left[g\left(X_0\right)\right] 
$$
The size of the error in the above limit decreases as the rate of decay
of the time autocorrelation
\begin{multline*}
\mathbf{corr}\left[g\left(Y_0^n\right),g\left(Y_0^0\right)\right] = \\
\frac{\mathbf{E}\left[
\left(g\left(Y_0^n\right)-\mathbf{E}\left[g\left(X_0\right)\right]\right)
\left(g\left(Y_0^0\right)-\mathbf{E}\left[g\left(X_0\right)\right]\right)\right]}
{\mathbf{Var}\left[g\left(X_0\right)\right]}
\end{multline*}
increases.
In this formula, $Y_0^0$ is assumed to be drawn from $\pi_0$.

It is well known that judicious elimination of variables by
renormalization can reduce long range spatial correlations (see e.g. \cite{bdfn92}).
The variables are
removed by averaging out their effects on the full distribution.
If the original density is $\pi(\hat{x},\tilde{x})$ and we wish to remove the
$\tilde{x}$ variables, the distribution of the remaining $\hat{x}$ 
variables is given by the marginal density (see \cite{c03,l02})
\begin{equation}\label{def:marginal}
\overline{\pi}\left(\hat{x}\right) = \int \pi\left(\hat{x},\tilde{x}\right)d\tilde{x}
\end{equation}
The full distribution  can be factored
as
$$
\pi(\hat{x},\tilde{x}) = \overline{\pi}(\hat{x})\pi(\tilde{x}\vert \hat{x})
$$
where $\pi(\tilde{x}\vert \hat{x})$ is the 
conditional density of $\tilde{x}$ given $\hat{x}$.
Because they exhibit shorter correlation
lengths, the marginal distributions are useful in the acceleration of Markov
chain Monte Carlo methods.

With this in mind we consider a collection
of lower dimensional Markov chains $Y_i^n\in\mathbb{R}^{d_i}$ which have stationary distributions
$\pi_i(x_i)$ where $d_0>\dots >d_i$.  For each $i\leq L$ 
let $T_i$ be the transition
probability density of $Y^n_i$, i.e.
$T_i(x_i\rightarrow y_i)$ is the probability density of
$\left\{Y_i^{n+1}=y_i\right\}$ given that $\left\{Y_i^n=x_i\right\}$.
The $\left\{\pi_i\right\}$ are approximate marginal 
distributions.  For example, 
divide the $x_i$ variables into two subsets,
$\hat{x}_i\in\mathbb{R}^{d_{i+1}}$ and $\tilde{x}_i\in\mathbb{R}^{d_i-d_{i+1}}$,
so that $x_i=\left(\hat{x}_i,\tilde{x}_i\right)$.  
The $\tilde{x}_i$ variables 
represent the variables of $x_i$ that are removed by marginalization, i.e.
$$
\pi_{i+1}\left(\hat{x}_i\right) \approx \int 
\pi_i\left(\hat{x}_i,\tilde{x}_i\right)\ d\tilde{x}_i.
$$

After arranging these chains in parallel we have the larger process 
$$Y^n = \left(Y^n_0,\dots, Y^n_L\right)
\in  \mathbb R^{d_0}
\times\dots\times\mathbb R^{d_L}.$$
The probability density of a move
to $\left\{Y^{n+1}=y\right\}$ given that $\left\{Y^n=x\right\}$
for $x,y\in\mathbb R^{d_0}
\times\dots\times\mathbb R^{d_L}$ is given by
\begin{equation}\label{def:T}
T(x\rightarrow y) = \prod_{i=0}^L T_i(x_i\rightarrow y_i).
\end{equation}
Since
\begin{equation*}
\int \left(T(x\rightarrow y) \prod_{i=0}^L\pi_i\left(x_i\right)\right)
\ dx_0\dots dx_L
= 
\prod_{i=0}^L\pi_i\left(y_i\right)
\end{equation*}
the stationary distribution of $Y^n$ is
\begin{equation*}
\Pi\left(x_0,\dots,x_L\right) = \pi_0\left(x_0\right)
\dots\pi_{L}\left(x_{L}\right).
\end{equation*}

The next step in the construction is to allow interactions between the chains
$\left\{Y_i^n\right\}$ and to thereby pass information from the rapidly 
equilibrating chains on the lower dimensional spaces (large $i$) down
to the chain on the original space ($i=0$).
This is accomplished by swap moves.
In a swap move between levels $i$ and $i+1$, we take a $d_{i+1}$ dimensional subset, 
$\hat{x}_i$, of the $x_i$ variables 
and exchange them with
the $x_{i+1}$ variables.  The remaining $d_i-d_{i+1}\ $ $\tilde{x}_i$ variables are resampled
from the conditional distribution $\pi_i\left(\tilde{x}_i\vert x_{i+1}\right)$.
For the full chain, this swap takes the form of a move from $\left\{Y^n=x\right\}$ 
to $\left\{Y^{n+1}=y\right\}$ where
$$
x = \left(\dots,\hat{x}_i,\tilde{x}_i,x_{i+1},\dots\right)
$$
and
$$
y = \left(\dots,x_{i+1},\tilde{y}_i,\hat{x}_i,\dots\right).
$$
The ellipses represent components of $Y^n$ that remain unchanged in the 
transition
and $\tilde{y}_i$ is drawn from $\pi_i\left(\tilde{x}_i\vert x_{i+1}\right)$.

If these swaps are undertaken unconditionally, the resulting chain with 
equilibrate rapidly, but will not, in general, preserve the product distribution
$\Pi$.  To remedy this we introduce the swap acceptance probability 
\begin{equation}\label{def:A_i}
A_i =
\min\biggl\lbrace 1,\ \frac{ \overline{\pi}_i(x_{i+1})\pi_{i+1}(\hat{x}_i)}
{\overline{\pi}_i(\hat{x}_i)\pi_{i+1}(x_{i+1})}\biggr\rbrace.
\end{equation}
In this formula $\overline{\pi}_i$ is the function on $\mathbb{R}^{d_{i+1}}$
resulting from marginalization of $\pi_i$ as in equation \ref{def:marginal}.
Given that $\left\{Y^n=x\right\}$, the probability density of 
$\left\{Y^{n+1}=y\right\}$, after the proposal and
either acceptance with probability $A_i$ or rejection
with probability $1-A_i$, of a swap move, is given by
\begin{multline*}
S_i\left(x\rightarrow y\right) =
 \left(1-A_i\right)\ 
 \delta_{\left\{y=x\right\}}\\
+ A_i\ \pi_i(\tilde{y}_i \vert x_{i+1})\ 
\delta_{\left\{\left(\hat{y}_i,y_{i+1}\right)
=\left(x_{i+1},\hat{x}_i\right)\right\}}
\prod_{j\notin\left\{i,i+1\right\}}\delta_{\left\{y_j=x_j\right\}}
\end{multline*}
for $x,y \in \mathbb R^{d_0}
\times\dots\times\mathbb R^{d_L}$. $\delta$ is the Dirac delta function.

We have the following lemma.
\begin{lemma}
The transition probabilities $S_i$ satisfy the
detailed balance condition for the measure $\Pi,$ 
i.e.
\begin{equation*}
\Pi(x)\ S_i\left(x\rightarrow y\right) = 
\Pi(y)\ S_i\left(y\rightarrow x\right)
\end{equation*}
where $x,y \in \mathbb R^{d_0}
\times\dots\times\mathbb R^{d_L}.$
\end{lemma}
The detailed balance condition stipulates that the probability of
observing a transition $x\rightarrow y$ is equal to that
of observing a transition $y\rightarrow x$ and guarantees that
the resulting Markov Chain preserves the distribution $\Pi$.
Therefore, under general conditions, averages over a trajectory of $\left\{Y^n\right\}$ 
will converge to averages over $\Pi$.  Since
$$
\pi_0(x_0) = \int \Pi(x_0,\dots,x_L)\ dx_1\dots dx_L
$$
we can calculate averages over $\pi_0$ by taking averages over the trajectories 
of the first $d_0$ components of $Y^n$.

\section{``Exact'' approximation of acceptance probability}
Notice that the formula \ref{def:A_i} for $A_i$ requires the evaluation of
$\overline{\pi}_i$ at the points $\hat{x}_i,x_{i+1}\in\mathbb{R}^{d_{i+1}}.$
While the approximation of 
$\overline{\pi}_i$ by functions on $\mathbb{R}^{d_{i+1}}$
is in general a very difficult 
problem, its evaluation at a single point is often not 
terribly demanding.  In fact, in many cases, including the examples in this
paper, the $\hat{x}_i$ variables can be chosen so that the remaining
$\tilde{x}_i$ variables are conditionally independent given $\hat{x}_i.$

Despite these mitigating factors, the requirement that we evaluate
$\overline{\pi}_i$ before we accept any swap
is a little onerous.  Fortunately, and somewhat surprisingly, 
this requirement is not necessary.
In fact, standard strategies for approximating the
point values of the marginals yield Markov chains that themselves preserve 
the target measure.  Thus even a poor estimate of the ratio appearing
in \ref{def:A_i} can give rise to a method that is exact in the sense that the
resulting Markov chain will asymptotically sample the target measure.

To illustrate this point, we consider the following 
example of a swap move.  Assume that the current position
of the chain is $\left\{Y^n=x\right\}$ where
$$
x = \left(\dots,\hat{x}_i,\tilde{x}_i,x_{i+1},\dots\right)
$$
The following steps will result in either $\left\{Y^{n+1}=x\right\}$
or $\left\{Y^{n+1}=y\right\}$ where
$$
y = \left(\dots,x_{i+1},\tilde{y}_i,\hat{x}_i,\dots\right)
$$
and $\tilde{y}_i\in\mathbb{R}^{d_i-d_{i+1}}$.
\begin{enumerate}
\item Let $v^0=\tilde{x}_i$  
  and let $v^j\in\mathbb{R}^{d_i-d_{i+1}}$ 
  for $j=1,\dots,M-1$ be independent
  samples from $p_i(\ \cdot\ \vert\hat{x}_i),$
  where $p_i(\ \cdot\ \vert\hat{x}_i)$ is a 
  reference density conditioned by $\hat{x}_i$. 
  For example, $p_i(\ \cdot\ \vert\hat{x}_i)$
  could be a Gaussian approximation of 
  $\pi_i(\tilde{x}_i\vert \hat{x}_i)$.
  How $p_i$ is chosen depends on the problem at hand 
  (see numerical examples below).  
  In general $p_i(\ \cdot\ \vert\hat{x}_i)$ should be easily evaluated and 
  independently sampled, and it should ``cover'' $\pi_i(\ \cdot\ \vert \hat{x}_i)$
  in the sense that areas of $\mathbb{R}^{d_i}$ 
  where $\pi_i(\ \cdot\ \vert \hat{x}_i)$
  is not negligible
  should be contained in areas where $p_i(\ \cdot\ \vert\hat{x}_i)$ is not negligible. 
\item
Let $u^j\in\mathbb{R}^{d_i-d_{i+1}}$ for 
$j=0,\dots,M-1$ be independent
random variables sampled from $p_i(\ \cdot\ \vert x_{i+1})$  (recall
that we are considering a swap of $\hat{x}_i$ and ${x_{i+1}}$ which
live in the same space).
Notice that the
$\left\{u^j\right\}$ variables depend on $x_{i+1}$
 while the $\left\{v^j\right\}$ variables depend on $\hat{x}_i$.
\item Define the weights
$$
w_v^j = \frac{\pi_i\left(\hat{x}_i,v^j\right)}
{p_i\left(v^j\vert\hat{x}_i\right)}
\ \ \text{and}\ \ 
w_u^j = \frac{\pi_i\left(x_{i+1},u^j\right)}
{p_i\left(u^j\vert x_{i+1}\right)}
$$
The choice of $p_i$ made above affects the variance of these weights, and therefore
the variance of the acceptance probability below.
\item Choose $\tilde{y}_i$ from among the $\left\{u^j\right\}$ according to the
multinomial distribution with probabilities
$$
\mathbf{P}\left(\tilde{y}_i=u^j\right) 
= \frac{w_u^j}{\sum_{l=0}^{M-1} w_u^l}.
$$
Notice that $\tilde{y}_i$ is an approximate sample from $\pi_i(\ \cdot\ \vert x_{i+1}).$
\item Set
$$Y^{n+1} = \left(\dots,x_{i+1},\tilde{y}_i,\hat{x}_i,\dots\right)$$
with probability
\begin{equation}\label{def:A^M_i}
A^M_i = \min\biggl\lbrace 1,\ \frac{ \pi_{i+1}(\hat{x}_i)\sum_{j=0}^{M-1} 
w^j_u}
{\pi_{i+1}(x_{i+1})
\sum_{j=0}^{M-1} w^j_v}\biggr\rbrace
\end{equation}
and
$$Y^{n+1} = Y^n=\left(\dots,\hat{x}_i,\tilde{x}_i,x_{i+1},\dots\right)$$
with probability
$
1-A^M_i
$.
\end{enumerate}

The transition probability density for the above swap move from
$x\rightarrow y$ for $x,y \in \mathbb R^{d_0}
\times\dots\times\mathbb R^{d_L}$is given by
\begin{multline*}
S^M_i(x\rightarrow y) = 
\left(1-R\right)\ 
 \delta_{\left\{y=x\right\}}\\
+R\  \delta_{\left\{\left(\hat{y}_i,y_{i+1}\right)
=\left(x_{i+1},\hat{x}_i\right)\right\}}
\prod_{j\notin\left\{i,i+1\right\}}\delta_{\left\{y_j=x_j\right\}}
\end{multline*}
where
\begin{multline*}
R = M \int p_i(u^0\vert x_{i+1}) 
\frac{w^0_u}{\sum_{j=0}^{M-1} w^j_u}\ A^M_i\\
\times\prod_{j=1}^{M-1} p_i(v^j\vert\hat{x}_i)p_i(u^j\vert x_{i+1})dv^j du^j
\end{multline*}
and $\delta$ is again the Dirac delta function.
In other words, $S^M_i$ dictates that the Markov chain accepts the swap with probability
$R$ and rejects it with probability $1-R$.

While the preceding swap move corresponds to a
method for approximating the ratio
$$
 \frac{ \overline{\pi}_i(x_{i+1})}
{\overline{\pi}_i(\hat{x}_i)}
$$
appearing in the formula for $A_i$ above, it also
has some similarities with
the multiple-try Metropolis method
presented in \cite{llw00} which uses multiple suggestion samples
to improve acceptance rates of standard MCMC methods.
The following lemma is suggested by results in \cite{llw00}.
\begin{lemma}
The transition probabilities $S^M_i$ satisfy the
detailed balance condition for the measure $\Pi.$ 
\end{lemma}
As before the detailed balance condition guarantees that averages over trajectories
of the first $d_0$ dimensions of $Y^n$ will converge to averages
over $\pi_0$.

The $A^M_i$ contain an approximation to the ratio 
of marginals in \ref{def:A_i}
\begin{align*}
\frac{\sum_{j=0}^{M-1} w_u^j}{\sum_{j=0}^{M-1} w_v^j}
=&\frac{\frac{1}{M}\sum_{j=0}^{M-1} \frac{\pi_i\left(x_{i+1},u^j\right)}
{p_i\left(u^j\vert x_{i+1}\right)}}
{\frac{1}{M}\sum_{j=0}^{M-1} \frac{\pi_i\left(\hat{x}_i,v^j\right)}
{p_i\left(v^j\vert \hat{x}_i\right)}}\\
\xrightarrow[{M\rightarrow\infty}]{a.s.}
&\frac{\mathbf{E}_{p_i}\left[\frac{\pi_i\left(x_{i+1},\widetilde{X}_i\right)}
{p_i\left(\widetilde{X}_i\vert x_{i+1}\right)}\ \vert
\left\{\widehat{X}_i=x_{i+1}\right\}\right]}
{\mathbf{E}_{p_i}\left[\frac{\pi_i\left(\hat{x}_i,\widetilde{X}_i\right)}
{p_i\left(\widetilde{X}_i\vert \hat{x}_i\right)}\ \vert
\left\{\widehat{X}_i=\hat{x}_i\right\}\right]}\\
=& \frac{\overline{\pi}_i(x_{i+1})}{\overline{\pi}_i(\hat{x}_i)}
\end{align*}
where $\mathbf{E}_{p_i}$ denotes expectation with respect to the density $p_i.$
When $0<\mathbf{E}_{p_i}\left[ w_v^j\ \vert
\left\{\widehat{X}_i=\hat{x}_i\right\}\right]<\infty$,
the convergence above follows from the strong law of large numbers
and the fact that
\begin{multline*}
\mathbf{E}_{p_i}\left[\frac{\pi_i\left(\widehat{X}_i\vert 
\widetilde{X}_i\right)}
{p_i\left(\widetilde{X}_i\vert \widehat{X}_i\right)}
\ \vert \left\{\widehat{X}_i=\hat{x}_i\right\}\right] 
= \int \frac{\pi_i(\hat{x}_i,\tilde{x}_i)}{p_i(\tilde{x}_i\vert\hat{x}_i)}
p_i(\tilde{x}_i\vert\hat{x}_i)\ d\tilde{x}_i\\
= \int \pi_i(\hat{x}_i,\tilde{x}_i)\ d\tilde{x}_i
= \overline{\pi_i}(\hat{x}_i)
\end{multline*}
For small values of $M$ in \ref{def:A^M_i}, calculation of the swap acceptance probabilities
is very cheap. However, higher values of $M$ 
may improve the acceptance rates.  For example, if the $\left\{\pi_i\right\}_{i>0}$ are
exact marginals of $\pi_0,$ then $A_i\equiv 1$ while $A^M_i\leq 1.$
Results similar to Lemma 2 hold when more general approximations replace the one
given above; for example when the $\left\{u^j\right\}$ and  $\left\{v^j\right\}$
are generated by a Metropolis-Hastings rule.  In practice one has to balance
the speed of evaluating  $A_i^M$ for small $M$ with the possible higher acceptance rates
for $M$ large.

It is easy to see that a Markov chain which evolves only by swap moves will only 
sample a finite number of configurations.  
These swap moves must therefore  
be used in conjunction with a transition rule that can
reach any region of space, such 
as $T$ from expression \ref{def:T}.  More precisely,
 $T$ should be $\Pi$-irreducible and aperiodic (see \cite{t94}).
The the transition rule for parallel marginalization is
\begin{multline*}
P(x\rightarrow y) = (1-\alpha)\ T(x\rightarrow y)\\
+ \alpha\ \int T(x\rightarrow z)S\left(z\rightarrow y\right)dz
\end{multline*}
where
$$
S(x\rightarrow y) = \sum_{k=0}^{L-1} \frac{1}{L} S^M_i\left(x\rightarrow y\right)
$$
and $\alpha\in \left[0,1\right)$ is the probability that a swap move occurs.
$P$ dictates that, with probability $\alpha$, the chain
attempts a swap move between levels $I$ and $I+1$ where $I$ is a random
variable chosen uniformly from $\left\{0,\dots,L-1\right\}$.
Next, each level of the chain evolves
independently according to the $\left\{T_i\right\}$.  
With probability
$1-\alpha$ the chain does not attempt a swap move, but does
evolve each level.
The next result follows trivially from Lemma 2 and guarantees the
invariance of $\Pi$ under evolution by $P$.
\begin{theorem}
The transition probability $P$ satisfies the
detailed balance condition for the measure $\Pi,$ 
i.e.
\begin{equation*}
\Pi(x)\ P\left(x\rightarrow y\right) = 
\Pi(y)\ P\left(y\rightarrow x\right)
\end{equation*}
where $x,y \in \mathbb R^{d_0}
\times\dots\times\mathbb R^{d_N}.$
\end{theorem}
Thus by combining standard MCMC steps on each component
governed by the transition probability $T$, with swap steps between
the components governed by $S$, an MCMC method results
which not only uses information from rapidly equilibrating
lower dimensional chains, but is also convergent.

\section{Numerical example 1:  bridge path sampling}
In the bridge path sampling problem
we wish to approximate conditional expectations of the form
$$
\mathbf{E}\left[g\left(Z^s\right)\ \vert
\{{Z^0=z^-}\},\{{Z^T=z^+}\}\right]
$$
where $s\in\left(0,T\right)$
 and $\left\{Z^t\right\}$ is the real valued processes given by the solution of the stochastic
differential equation
\begin{equation}\label{sde1}
dZ^t = f\left(Z^t\right)dt
+ \sigma\left(Z^t\right)dW^t.
\end{equation}
$g$, $f$ and $\sigma$
 are real valued functions of $\mathbb{R}$.  Of course we can also consider functions
 $g$ of more than one time.
This problem arises, for example, in financial volatility estimation.
Because in general we cannot sample paths of \ref{sde1} we must
first approximate $\left\{Z^t\right\}$ by a discrete process 
for which the path density is readily available.
Let $t_0=0,t_1=\frac{T}{K},\dots,t_K=T$ be a mesh on which we wish to calculate path
averages.
One such approximate process is given by the 
linearly implicit
Euler scheme (a balanced implicit method, see \cite{mps98}),
\begin{equation}\label{def:lie}
\begin{split}
&X^{t_{k+1}} = X^{t_k} + f\left(X^{t_k}\right)\triangle\\
&\ \ \ \ + \left(X^{t_{k+1}} - X^{t_k}\right)f^{'}\left(X^{t_k}\right)\triangle
+ \sigma\left(X^{t_k}\right)\sqrt{\triangle}\ \xi^k,\\
&X^0 = Z^0\ \ \ \ X^{t_K}=Z^T.
\end{split}
\end{equation}
The $\left\{\xi^k\right\}$ are independent Gaussian random
variables with mean 0 and variance 1,
and $\triangle=\frac{T}{K}.$  $K$ is assumed to be a power
of 2.
The choice of this scheme over the Euler scheme (see \cite{kp92}) 
is due to its favorable stability properties as explained later.
Without the condition $X^{t_K}=Z^T$ above, generating
samples of 
$\left(X^0,\dots,X^{t_K}\right)$ is a relatively straitforward
endeavor.
One simply generates a sample of $Z^0$,
then evolves the system with this initial condition.  
However, the presence of
information about $\left\{Z^t\right\}_{t>0}$ complicates 
the task.
In general, some sampling method which requires only 
knowlege of a function
proportional to conditional density of  
$\left(X^{t_1},\dots,X^{t_{K-1}}\right)$
must be applied.
The approximate path density associated with discretization
\ref{def:lie} is
\begin{multline}\label{def:pi_0:ex1}
\pi_0\left(x^{t_1},\dots,x^{t_{K-1}}\ \vert x^0,x^{t_K}\right) \propto\\
\exp\left(-\sum_{k=0}^{K-1} 
V\left(x^{t_k},x^{t_{k+1}},\triangle\right)
\right)
\end{multline}
where
\begin{equation*}\label{def:V}
V\left(x,y\right) = \\
\frac{\left[\left(1-\triangle f^{'}\left(x\right)\right)\left(y
-x\right)+\triangle f\left(x\right)\right]^2}
{2\sigma^2\left(x\right)\triangle}
\end{equation*}

At this point we wish to apply the parallel marginalization sampling
procedure to the density $\pi_0$.
However, as indicated above, a prerequisite for the use of parallel 
marginalization is the ability to estimate marginal densities.
In some important problems
 homogeneities in the underlying system yield 
simplifications in the calculation of these densities by the
methods in \cite{c03,s05}.
These calculations
can be carried out before implementation of parallel marginalization,
or they can be integrated into the sampling procedure.

In some cases, the numerical estimation of the $\left\{\pi_i\right\}_{i>0}$
can be completely avoided.  
The examples presented here are two such cases.
Let $S_i = \left\{0,2^i,\dots,K\right\}$.  Decompose $S_i$
as $\widehat{S}_i\sqcup\widetilde{S}_i$ where
$$\widehat{S}_i = \left\{0,2(2^i),4(2^i),\dots,K\right\}$$
 and
$$\widetilde{S}_i = \left\{2^i,3(2^i),5(2^i),\dots,K-2^i\right\}.$$
In the notation of the previous sections, $x_i=\left(\hat{x}_i,\tilde{x}_i\right)$
where $\hat{x}_i=\left\{x_i^{t_k}\right\}_{k\in\widehat{S}_i\setminus\left\{0,K\right\}}$ and
$\tilde{x}_i=\left\{x_i^{t_k}\right\}_{k\in\widetilde{S}_i}.$
In words, the hat and tilde variables represent alternating time slices of the path.
For all $i$ fix $x^0_i=z^-$ and $x_i^{t_K}=z^+$.
We choose the approximate marginal densities
$$
\pi_i\left(\left\{x_i^{t_k}\right\}_{k\in S_i\setminus\left\{0,K\right\}}\ 
\vert x_i^0,x_i^{t_K}\right) 
\propto q_i\left(\left\{x_i^{t_k}\right\}_{k\in S_i}\right)
$$ 
where
for each $i$, $q_i$ is defined by
successive coarsenings of \ref{def:lie}.
That is,
\begin{multline*}\label{def:q_i}
q_i\left( \left\{x_i^{t_k}\right\}_{k\in S_i} \right)\\
= \exp\left(-\sum_{k=0}^{K/2^i-1} 
V\left(x_i^{t_{2^ik}},x_i^{t_{2^i(k+1)}},2^i\triangle\right)\right).
\end{multline*} 
Since $\pi_i$ will be sampled
using a Metropolis-Hastings method 
with  $x^0$ and $x^{t_K}$ fixed,
knowlege of the normalization constants
$$
\mathcal{Z}_i\left(x_i^0,x_i^{t_K}\right) 
= \int q_i\ \prod_{k\in S_i\setminus\left\{0,K\right\}}dx_i^{t_k}
$$
is unnecessary.

Notice from \ref{def:pi_0:ex1} that, conditioned on the values of $x^{t_{k-1}}$ and
$x^{t_{k+1}}$, the variance of $x^{t_k}$ is of order $\triangle$.  Thus
any perturbation of $x^{t_k}$ which leaves $x^{t_j}$ fixed for $j\neq k$ and
which is compatible with joint distribution \ref{def:pi_0:ex1} must be of the
order $\sqrt{\triangle}$.  This suggests that distributions defined by
coarser discretizations of \ref{def:pi_0:ex1} will allow larger perturbations,
and consequently will be easier to sample.  However, it is important to 
choose a discretization that remains stable for large values of $\triangle$.
For example, while the linearly implicit Euler method performs well in 
the experiments below, similar tests
using the Euler method were less
successful due to limitations on the largest allowable values of $\triangle$.

In this numerical example 
bridge paths are sampled between time 0 and time 10 for a diffusion in 
a double well potential
$$
f(x) = -4x\left(x^2-1\right)\ \ \text{and}\ \ \sigma(x) = 1
$$
The left and right end points are chosen
as $z^-=z^+=0$.
$Y^n_i\in\mathbb{R}^{10/\left(2^i\triangle\right)+1}$
is the $i^{th}$ level of the parallel marginalization 
Markov chain at algorithmic time $n$.
There are 10 chains ($L = 9$ in expression \ref{def:T}).
The observed swap acceptance rates 
are reported in Table 1.
Let $Y^n_{mid}\in\mathbb{R}$ denote the midpoint of the path defined by $Y^n_0$
(i.e. an approximate sample of the path at time 5).  
In Fig. 1 the autocorrelation of $Y^n_{mid}$
$$
\mathbf{corr}\left[Y^n_{mid},Y^0_{mid}\right]
$$
is compared to that of a standard Metropolis-Hastings rule.
In the figure, the time scale of the autocorrelation for the Metropolis-Hastings method has
been scaled by a factor of 1/10 to more than account for the extra computational time
required per iteration of parallel marginalization.
The relaxation time of the parallel chain is clearly reduced.
In these numerical examples, the algorithm in the previous section is applied
with a slight simplification. 
First generate M independent Gaussian random paths 
$\left\{\zeta^j\left(t_k\right)\right\}_{k\in\widetilde{S}_i}$  with independent components
$\zeta^j\left(t_k\right)$ of mean 0 and variance $2^{i-1}\triangle$.
 For each $j$ and $k\in\widetilde{S}_i$ let 
$$
u^j\left(t_k\right) = \zeta^j\left(t_k\right)+0.5\left(x^{t_{k-1}}_{i+1}+x^{t_{k+1}}_{i+1}\right)
$$
If in step 4, $\tilde{y}_i=u^{j^*}$, then in step 1 we set $v^0=\tilde{x}_i$ and 
for each $k\in\widetilde{S}_i$
$$\left\{v^j\left(t_k\right)\right\}_{j>0}=
\left\{\zeta^j\left(t_k\right)+0.5\left(\hat{x}_i^{t_{k-1}}+\hat{x}_i^{t_{k+1}}\right)\right\}_{j\neq j^*}.
$$
All other steps remain the same.  This change yields a slightly faster though less generally 
applicable swap step that also preserves the density $\Pi$.
Notice that this modification implies that the reference density $p_i$ is given by
$$
p_i\left(\tilde{x}_i\vert \hat{x}_i\right)
\propto
\exp\left(
\sum_{k\in\widetilde{S}_i} 
-\frac{\left(
\tilde{x}^{t_k}_i
-0.5\left(\hat{x}_i^{t_{k-1}}+\hat{x}_i^{t_{k+1}}\right)
\right)^2}
{2^i\triangle}
\right).
$$
For this
problem, the choice of $M$ in \ref{def:A^M_i}, 
the number of samples of $\left\{u^j\right\}$ 
and $\left\{v^j\right\}$, seems to have little effect on the
swap acceptance rates.  In the numerical experiment $M=i+1$
for swaps between levels $i$ and $i+1$.

\section{Numerical example 2: non-linear smoothing/filtering}
In the non-linear smoothing and filtering problem
we wish to approximate conditional expectations of the form
$$
\mathbf{E}\left[g\left(Z^s\right)  \ \vert
\left\{H^j=h^j\right\}_0^J\right]
$$
where $s\in\left(0,T\right)$
 and the real valued processes
$\left\{Z^t\right\}$ and $\left\{H^j\right\}$ are given by
the system
\begin{equation*}
\begin{split}
&dZ^t = f\left(Z^t\right)dt
+ \sigma\left(Z^t\right)dW^t, \\
&H^j = r\left(Z^{s_j}\right) + \chi^j,\\
&Z^0 \sim \rho,\ \ \ \ \chi^n\sim i.i.d.\ \mu.
\end{split}
\end{equation*}
$g$, $f$, $\sigma$, and $r$ are real valued functions of $\mathbb{R}$.
The $\left\{\chi^j\right\}$ are
real valued independent random variable drawn from the density $\mu$ and are
independent of the Brownian motion $\left\{W^t\right\}.$
$\left\{s_j\right\}\subset\left\{t_j\right\},$ 
and $0=s_0<s_1<...<s_J=T.$
The process $Z^t$ is a hidden signal and the $\left\{H^j\right\}$ are 
noisy observations.

Again, the system must first be discretized.  The
linearly implicit Euler scheme gives
\begin{equation*}\label{sde2}
\begin{split}
&X^{t_{k+1}} = X^{t_k} + f\left(X^{t_k}\right)\triangle\\
&\ \ \ \ + \left(X^{t_{k+1}} - X^{t_k}\right)f^{'}\left(X^{t_k}\right)\triangle
+ \sigma\left(X^{t_k}\right)\sqrt{\triangle}\ \xi^k,\\\
&H^j = r\left(X^{s_j}\right)+\chi^j,\notag\\
&X^0 = Z^0\ \ \ \ \chi^n\sim i.i.d.\ \mu.
\end{split}
\end{equation*}
The $\left\{\xi^k\right\}$ are independent Gaussian random
variables with mean 0 and variance 1,
and $\triangle=\frac{T}{K}.$
  The $\left\{\xi^k\right\}$ are independent of the
$\left\{\chi^j\right\}$. $K$ is again assumed to be a power of 2.

The approximate path measure for this problem is
\begin{multline*}
\pi_0 \left(x^{t_0},\dots,x^{t_K}\ \vert
h^0,\dots,h^T\right)\propto\\
\exp\left(
-\sum_{k=0}^{K-1}
V\left(x^{t_k},x^{t_{k+1}},\triangle\right)
\right)\\
\times\rho\left(x_i^{t_0}\right)\prod_{n=0}^{J} 
\mu\left(x_i^{s_j}-r\left(h^{s_j}\right)\right)
\end{multline*}
The approximate marginals are chosen as
\begin{multline*}
\pi_i \left(\left\{x_i^{t_k}\right\}_{k\in S_i}\ 
\vert h^0,\dots,h^T\right)
\propto\\
q_i\left(\left\{x_i^{t_k}\right\}_{k\in S_i}\right)
\rho\left(x^{t_0}\right)\prod_{n=0}^{J} \mu\left(x^{s_j}-r\left(h^{s_j}\right)\right)
\end{multline*}
where $V$, $q_i$ and $S_i$ are as defined in the previous section.

In this example, samples of the smoothed path are generated
between time
time 0 and time 10 for the same diffusion in 
a double well potential.  The densities $\mu$ and $\rho$ are chosen 
as
$$
\mu = N(0,0.01)\ \ \text{and}\ \ 
\rho(x) \propto \exp\left(-\left(x^2-1\right)^2\right)
$$
The observation times are $s_0 = 0,s_1=1,\dots,s_{10}=10$
with
$H^j=-1$ for $j=0,\dots,5$ and $H^j=1$ for $j=6,\dots,10$.
There are 8 chains ($L = 7$ in expression \ref{def:T}).
The observed swap acceptance rates 
are reported in Table 1.
Again, $Y^n_{mid}\in\mathbb{R}$ denotes the
 midpoint of the path defined by $Y^n_0$
(i.e. an approximate sample of the path at time 5).  
In Fig. 2 the autocorrelation of $Y^n_{mid}$
is compared to that of a standard Metropolis-Hastings rule.
The figure has been adjusted as in the previous example.
The relaxation time of the parallel chain is again 
clearly reduced. The algorithm is modified as in the previous example. 
For this problem, acceptable 
swap rates require a higher choice of $M$ in \ref{def:A^M_i} than needed
in the bridge sampling problem.  In this numerical experiment $M=2^i$
for swaps between levels $i$ and $i+1$. 

\section{Conclusion}
A Markov chain Monte Carlo method has been proposed and applied to two 
conditional path sampling problems for stochastic differential equations.
Numerical results indicate that this method, parallel marginalization,
can have a dramatically reduced equilibration time when compared to standard MCMC 
methods.  

Note that parallel marginalization should not be viewed as 
a stand alone method.  Other acceleration techniques such as hybrid Monte Carlo 
can and should be implemented at each level within
the parallel marginalization framework.
As the smoothing problem indicates, the acceptance probabilities at
coarser levels can become small.  The remedy for this is the development of 
more accurate approximate marginal distributions by, for example,
the methods
in \cite{c03} and \cite{s05}.

\begin{acknowledgments}
I would
like to thank Prof. A. Chorin for his guidance during this research, which
was carried out during my Ph.D. studies at U. C. Berkeley.
I would also like to thank Dr. P. Okunev, and Dr. P. Stinis
for very helpful discussions and comments.  
This work was supported by the Director, Office of Science, Office
of Advanced Scientific Computing Research, of the U. S. Department of
Energy under Contract No. DE-AC03-76SF00098.
\end{acknowledgments}




\end{article}



\begin{figure}[t]
\begin{center}
\includegraphics[width=.4\textwidth]{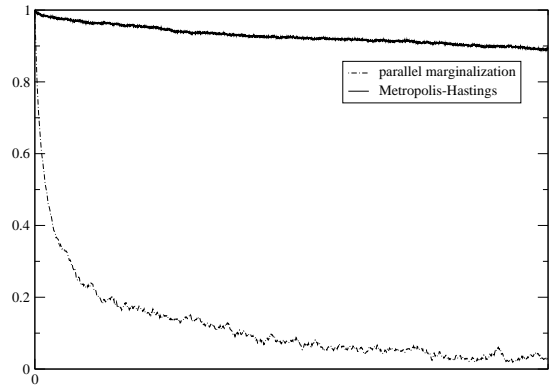}
\caption{Autocorrelations of parallel marginalization and standard
Metropolis-Hastings methods for bridge sampling problem.}\label{xafoto}
\end{center}
\end{figure}

\begin{figure}[t]
\begin{center}
\includegraphics[width=.4\textwidth]{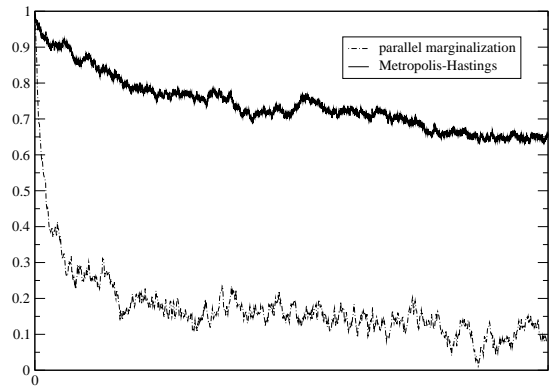}
\caption{Autocorrelations of parallel marginalization and standard
Metropolis-Hastings methods for filtering and smoothing
problem.}\label{xafoto}
\end{center}
\end{figure}

\begin{table*}[b]
\caption{Swap acceptance rates for bridge sampling and filtering/smoothing
problems}
\begin{tabular*}{\hsize}{@{\extracolsep{\fill}}rrrrrrrrrrrrr}
\hline
\multicolumn1c{Levels\tablenote{Swaps between levels $i$ and $i+1$}}&
\multicolumn1c{0/1}&\multicolumn1c{1/2}&
\multicolumn1c{2/3}&\multicolumn1c{3/4}&
\multicolumn1c{4/5}&\multicolumn1c{5/6}&
\multicolumn1c{6/7}&\multicolumn1c{7/8}
&\multicolumn1c{8/9}
\cr\hline
\multicolumn1c{BS\tablenote{Bridge sampling problem}}&
\multicolumn1c{0.86}&\multicolumn1c{0.83}&
\multicolumn1c{0.75}&\multicolumn1c{0.69}&
\multicolumn1c{0.54}&\multicolumn1c{0.45}&
\multicolumn1c{0.30}&\multicolumn1c{0.22}
&\multicolumn1c{0.26}
\cr
\multicolumn1c{FS\tablenote{Filtering/smoothing problem}}&
\multicolumn1c{0.86}&\multicolumn1c{0.83}&
\multicolumn1c{0.74}&\multicolumn1c{0.65}&
\multicolumn1c{0.46}&\multicolumn1c{0.23}&
\multicolumn1c{0.04}&\multicolumn1c{NA}
&\multicolumn1c{NA}
\cr
\hline
\end{tabular*}
\tablenotes
\end{table*}






\end{document}